\documentclass[]{amsart}
\usepackage[]{amsmath, amsthm, amsfonts, amssymb}
\usepackage[]{graphicx}
 \usepackage[all]{xy}

 \newcommand {\C} {{\mathbb C}}
 
 \newcommand {\Z} {{\mathbb Z}}
 \newcommand {\Q} {{\mathbb Q}}
\newcommand{\dt} {{\bullet}}
\newcommand {\dd} {{\partial}}
\newcommand {\db} {{\bar \partial}}
\newcommand {\F} {{\mathcal F}}
\newcommand {\OO} {{\mathcal O}}
\newcommand {\X}  {{\mathcal X}}

\newtheorem{thm}[subsection]{Theorem}
 \newtheorem{cor}[subsection]{Corollary}
 \newtheorem{lemma}[subsection]{Lemma}
 \newtheorem{prop}[subsection]{Proposition}
 
 \newtheorem{rmk}[subsection]{Remark}
 \newtheorem{ex}[subsection]{Example}

\begin{document}
\title{The combinatorial part of the cohomology  of a singular variety }
\author{Donu Arapura}
\author{Parsa Bakhtary}
\author{Jaros{\l}aw W{\l}odarczyk}
\thanks{First and third author partially supported by NSF}
 \address{Department of Mathematics\\
   Purdue University\\
   West Lafayette, IN 47907\\
   U.S.A.}
\date{}
\maketitle

\begin{abstract}
  We study the first step of the weight filtration on the cohomology
  of a proper complex algebraic variety, which we call the combinatorial part. 
  We obtain a natural upper  bound on its size, which gives rather strong information about the topology of rational singularities.
\end{abstract}

Given a possibly reducible complex algebraic variety $X$, we define the
{\em combinatorial part} of the compactly supported cohomology
to  a subspace $KH_c^i(X)\subseteq H_c^i(X,\Z)$ characterized by the following
axioms:
\begin{enumerate}
\item[(K1)] These subspaces are preserved by proper pullbacks.
\item[(K2)] If $X$ is smooth and complete, $KH_c^0(X) = H_c^0(X)$
and $KH_c^i(X) = 0$ for $i>0$.
\item[(K3)] If  $U\subseteq X$ is an open immersion and $Z=X-U$, then
the standard exact sequence
$$\ldots H_c^{i-1}(Z)\to H_c^i(U)\to H_c^i(X)\to \ldots$$
restricts to an exact sequence
$$\ldots KH_c^{i-1}(Z)\to KH_c^i(U)\to KH_c^i(X)\to \ldots$$
\end{enumerate}

The proof of uniqueness, when $X$ is complete, given below is a simple induction. Existence will follow
by identifying $KH_c^i(X)$  with the 
first step of the weight filtration $W_0H_c^i(X)$ of Deligne \cite{deligne} and Gillet-Soul\'e \cite{gs}.
It will be both convenient and necessary
to review the basic construction which gives a method for calculating
this in terms of the underlying combinatorics of a simplicial resolution.
In simple cases, such as when $X$ has simple normal crossing singularities,
this can be made quite explicit. We note that in this paper varieties are
reduced schemes of finite type over $\C$.
We can extend this an arbitrary complex scheme of finite type $X$,
by defining $KH_c^i(X) = KH_c^i(X_{red})$.

 Work of Stepanov \cite{step} and the second author
\cite{parsa} suggested a certain natural bound on the dimension of the
combinatorial part of cohomology of the exceptional divisor of a
singularity. The main purpose of this note is to verify this in a refined
form. Given a proper map of varieties $f:X\to Y$, we show that $\dim
KH^i(f^{-1}(y))$ is bounded above by $\dim (R^if_*\OO_X)_y\otimes
\OO_y/m_y$. In particular, in accordance with a conjecture of
Stepanov, the first space vanishes for a resolution of a rational
singularity.

\section{Uniqueness for complete varieties}

As a warm up, we prove the uniqueness statement for complete varieties.
For this it is convenient to replace (K3) by (K3$'$) below.

\begin{lemma}\label{lemma:mayer}
 Assume that $KH_c^i$ satisfies the  axioms (K1)-(K3).
 Given a   complete variety $X$
with closed set $S$ and a desingularization $f:\tilde X\to X$
which is an isomorphism over $X-S$. Let $E=f^{-1}(S)$.
\begin{enumerate}
\item[(K3$'$)]
 Then there is an exact sequence
$$
\ldots \to KH^{i-1}(E)\to KH^i(X)\to
KH^{i}(\tilde X)\oplus KH^{i}(S)\to \ldots
$$

\end{enumerate}

\end{lemma}

\begin{proof}
 This   follows from  diagram chase on
$$
\xymatrix{
\ldots KH^{i-1}(S)\ar[r]\ar[d] & KH_c^i(U)\ar[r]\ar[d]^{=} & KH^i(X)\ar[r]\ar[d] & KH^i(S)\ar[d]\ldots \\ 
 KH^{i-1}(E)\ar[r] & KH_c^i(U)\ar[r] & KH^i(\tilde X)\ar[r] & KH^i(E)
}
$$
\end{proof}

\begin{rmk}\label{rmk:mayer}
  In general, by the same argument we get a sequence
$$
\ldots \to KH_c^{i-1}(E)\to KH_c^i(X)\to
KH_c^{i}(\tilde X)\oplus KH_c^{i}(S)\to \ldots
$$
\end{rmk}

\begin{lemma}\label{lemma:mayer2}
 Assume that $KH^i$ satisfies the  axioms (K1), (K2) and (K3$'$).
 Given a complete variety $X$
with  a closed set $S$ and a desingularization $f:\tilde X\to X$
which is an isomorphism over $X-S$. Let $\tilde S\to S$ be a desingularization
of $S$ and $F= \tilde S\times_X \tilde X$. Then
there is an exact sequence
$$
\ldots \to KH^{i-1}(F)\to KH^i(X)\to
KH^{i}(\tilde X)\oplus KH^{i}(\tilde S)\to \ldots
$$
\end{lemma}

\begin{proof}
  Consider the diagram
$$
\xymatrix{
 \tilde F\ar[rd]^{f}\ar[dd]\ar[rr] &  & \tilde X\times  \tilde S\ar[rd]^{p}\ar[dd]&  \\ 
  & F\ar[dd]\ar[rr] &  & \tilde X\ar[dd] \\ 
 \tilde S\ar[rd]^{id_{\tilde S}}\ar[rr]^{\gamma} &  & X\times  \tilde S\ar[rd]^{p} &  \\ 
  & \tilde S\ar[rr] &  & X
}
$$
where the maps labelled $p$ are projections, $\gamma$ is the graph
of the composition $\tilde S\to S\to X$, and 
$\tilde F= \tilde S\times_{X\times \tilde  S} (\tilde X\times  \tilde S)$.
The lefthand square containing $f$ and $id_s$ is easily seen to be Cartesian.
Therefore $f$ gives an isomorphism
$\tilde F\cong F$.  Thus from the previous lemma, we obtain an 
exact sequence
$$
\ldots \to KH^{i-1}(F)\to KH^i(X\times  S)\to
KH^{i}(\tilde X\times  S)\oplus KH^{i}(\tilde S)\to \ldots
$$
Choose base points $s_1,\ldots s_N$ in each connected component of $ S$,
Define $\sigma=\frac{1}{N}\sum_j (id\times s_j)*: H^i(X\times S)\to H^i(X)$.
This gives a left inverse to $p^*$.  Then a  diagram chase using 
$$
\xymatrix{
 KH^{i-1}(F)\ar@{-}[d]^{=}\ar[r]^{\delta} & KH^i(X\times  S)\ar[d]^{\sigma}\ar[r] & KH^{i}(\tilde X\times  S)\oplus KH^{i}(\tilde S)\ar[r]\ar[d] & KH^{i}(F)\ar@{-}[d]^{=} \\ 
 KH^{i-1}(F)\ar[r]^{\sigma\circ\delta} & KH^i(X)\ar@<1ex>[u]^{p^*}\ar[r] & KH^{i}(\tilde X)\oplus KH^{i}(\tilde S)\ar@<1ex>[u]\ar[r] & KH^{i}(F)
}
$$
shows that the bottom row is exact.

\end{proof}

\begin{thm}
There is at most one collection of subspaces   $KH^i(X)\subseteq H^i(X,\Z)$,
with $X$ complete, 
satisfying  axioms (K1), (K2) and (K3$'$).
\end{thm}

\begin{proof}
We prove this by induction on $i$.
  First we check that  $KH^0(X) = H^0(X)$. We can assume that $X$ is 
connected. 
If $p\in X$, then 
$$H^0(X) = H^0(p)=KH^0(p)= KH^0(X)$$
by the axioms. 

By lemma~\ref{lemma:mayer2}, there is an exact sequence
$$
KH^{i-1}(F)\to KH^i(X)\to
KH^{i}(\tilde X)\oplus KH^{i}(\tilde S)=0
$$
So $KH^i(X) = im[KH^{i-1}(F)\to H^i(X)]$.
\end{proof}

\section{Simplicial resolutions}

The general construction is based on simplicial resolutions.
We start by recalling some standard material \cite{deligne, gnpp, ps}.
A simplicial object in a category is a diagram
$$
\xymatrix{ \ldots X_2\ar[r]\ar@<1ex>[r]\ar@<-1ex>[r] &
  X_1\ar@<1ex>[r]\ar@<-1ex>[r] & X_0 }
$$
with $n$ face maps $\delta_i:X_n\to X_{n-1}$ satisfying the standard
relation $\delta_i\delta_j = \delta_{j-1}\delta_i$ for $i<j$; this
would be more accurately called a ``strict simplicial'' or
``semisimplicial'' object since we do not insist on degeneracy maps
going backwards.  The basic example of a simplicial set,
i.e. simplicial object in the category of sets, is given by taking
$X_n$ to be the set of $n$-simplices of a simplicial complex on an
ordered set of vertices. Let $\Delta^n$ be the standard $n$-simplex
with faces $\delta_i':\Delta^{n-1}\to \Delta^n$.  Given a simplicial
set or more generally a simplicial topological space, we can glue the
$X_n\times \Delta^n$ together by identifying $(\delta_i x,y)\sim
(x,\delta_i' y)$. This leads to a topological space $|X_\dt|$ called
the geometric realization, which generalizes the usual construction of
the topological space associated a simplicial complex.

Given a simplicial space, filtering $|X_\dt|$ by skeleta
$\bigcup_{n\le N} X_n\times \Delta^n/\sim$ yields the spectral
sequence
\begin{equation}
  \label{eq:ss1}
  E_1^{pq} = H^q(X_p,A)\Rightarrow H^{p+q}(|X_\dt|,A)  
\end{equation}
for any abelian group $A$.
It is convenient to extend this. A simplicial sheaf on $X_\dt$ is a
collection of sheaves $\F_n$ on $X_n$ with ``coface'' maps
$\delta_i^{*}\F_{n-1}\to \F_n$ satisfying the face relations.  For
example, the constant sheaves $\Z_{X_\dt}$ with identities for coface
maps forms a simplicial sheaf.  If $X_\dt$ is a simplicial object in
the category of complex manifolds, then $\Omega^i_{X_\dt}$ with the
obvious maps, forms a simplicial sheaf. We can define cohomology by
setting
$$H^i(X_\dt,\F_\dt) = Ext^i(\Z_{X_\dt},\F_\dt)$$
This generalizes sheaf cohomology in the usual sense, and it can be
extended to the case where $\F_\dt^\dt$ is a bounded below complex of
simplicial sheaves by using a hyper $Ext$.  When $\F=A$ is constant,
this coincides with $H^{i}(|X_\dt|,A)$. But in general the meaning is
more elusive.  There is a spectral sequence
\begin{equation}
  \label{eq:ss2}
  E_1^{pq}(\F^\dt_\dt) = H^q(X_p,\F^\dt_p)\Rightarrow H^{p+q}(X_\dt, \F^\dt_\dt)  
\end{equation}
generalizing \eqref{eq:ss1}. Filtering $\F^\dt$ by the ``stupid
filtration'' $\F_\dt^{\ge n}$ yields a different spectral sequence
\begin{equation}
  \label{eq:ss3}
  {}'E_1^{pq} = H^q(X_\dt,\F^p_\dt)\Rightarrow H^{p+q}(X_\dt, \F^\dt_\dt)  
\end{equation}

\begin{thm}[Deligne]
  If $X_\dt$ is a simplicial object in the category of compact
  K\"ahler manifolds and holomorphic maps. The spectral sequence
  \eqref{eq:ss1} degenerates at $E_2$ when $A=\Q$. 
\end{thm}

\begin{rmk}
  The theorem follows from a more general result in
  \cite[8.1.9]{deligne}. However the argument is very complicated.
  Fortunately, as pointed out in \cite{dgms}, this special case
  follows easily from the $\dd\db$-lemma.  Here we give a more
  complete argument.
\end{rmk}

\begin{proof}
  It is enough to prove this after tensoring with $\C$.  We can
  realize the spectral sequence as coming from the double
  $(E^\dt(X_\dt),d,\pm\delta)$, where $(E^\dt,d)$ is the $C^\infty$ de
  Rham complex, and $\delta$ is the combinatorial differential. (We
  are mostly going to ignore sign issues since they are not relevant
  here.) In fact this is a triple complex, since each $E^\dt(-)$ is
  the total complex of the double complex $(E^{\dt\dt}(-),\dd,\db)$.

  Given a class $[\alpha]\in H^i(X_j)$ lying in the kernel of
  $\delta$, we have $\delta\alpha = d\beta$ for some $\beta\in
  E^{i-1}(X_{j+1})$ Then $d_2([\alpha])$ is represented by
  $\delta\beta\in E^{i-1}(X_{j+2})$.  We will show this vanishes in
  cohomology. The ambiguity in the choice of $\beta$ will turn out to be 
  the key point.

  By the Hodge decomposition, we can assume that $\alpha$ is pure of
  type $(p,q)$. Therefore $\delta\alpha$ is also pure of this type. We
  can now apply the $\dd\db$-lemma \cite[p 149]{gh} to write $\alpha = \dd\db \gamma$
  where $\gamma\in E^{p-1,q-1}(X_{j+1})$. This means we have two
  choices for $\beta$.  Taking $\beta=\db\gamma$ shows that
  $d_2([\alpha])$ is represented by a form of pure type $(p-1,q)$. On
  the other hand, taking $\beta = -\dd\gamma$ shows that this class is
  of type $(p,q-1)$.  Thus $d_2([\alpha]) \in H^{p-1,q}\cap
  H^{p,q-1}=0$.

  By what we just proved $\delta\alpha = d\beta, \delta\beta=d\eta$,
  and $\delta\eta$ represents $d_3([\alpha])$. It should be clear that
  one can kill this and higher differentials in the exact same way.
 
\end{proof}

\begin{cor}
  With the same assumptions as the theorem, the spectral sequence
  \eqref{eq:ss2} degenerates at $E_2$ when $\F=\OO_{X_\dt}$.
\end{cor}

(This fixes an incorrect proof in \cite[2.4]{step}.)

\begin{proof}
  By the Hodge theorem, the spectral sequence for $\F=\OO_{X_\dt}$ is
  a direct summand of the spectral sequence for $\F = \C$.
\end{proof}

\begin{thm}[Deligne]
 Given any (possibly reducible) variety $X$, there
exists a smooth simplicial variety $X_\dt$, which 
we call a {\em simplicial resolution},
with proper morphisms $\pi_\dt:X_\dt\to X$ (commuting with face maps) inducing
a homotopy equivalence between $|X_\dt|$ and $X$. Given a morphism $f:X\to Y$
there exists simplicial resolutions $X_\dt, Y_\dt$ and a 
morphism $f_\dt:X_\dt\to Y_\dt$ compatible with $f$.
\end{thm}

 The theorem is a consequence of resolution of singularities. Proofs can be found 
in \cite{deligne,gnpp, ps}.  Note that the original construction of
Deligne results in a  necessarily infinite diagram, 
whereas the method of Guillen et. al
yields a fairly economical resolution.
Here are some examples.
 
\begin{ex}\label{ex:ncd}
  Suppose that $X$ is an analytic space, whose irreducible components
  $X^i$ are compact K\"ahler, and suppose that their intersections
  $X^{ij\ldots}=X^i\cap X^j\ldots$ are all smooth.  This includes 
  the case of a divisor with simple normal crossings. Then an explicit
  simplicial resolution is given by taking $X_i$ to be the disjoint
  union of $(i+1)$-fold intersection of components of $X$. The face
  map $\delta_k$ is given by inclusions
$$X^{i_1\ldots i_n}\subset X^{i_i\ldots\hat i_k\ldots i_n}\>( i_1<\ldots< i_k)$$ 
\end{ex}

\begin{rmk}\label{rmk:ncd}
  The above construction makes perfect sense for general $X$, and it yields a 
  (generally singular) simplicial variety $X_\dt$ with $|X_\dt|$ homotopic to $X$.  
  $X_\dt$ will always be dominated by a simplicial resolution.
\end{rmk}

\begin{ex}\label{ex:isolated}
  Following the method of \cite{gnpp} we can construct a simplicial
  resolution of a variety $X$ with isolated singularities as
  follows. Let $f:\tilde X\to X$ be a resolution of singularities such
  that the exceptional divisor $E = \cup E^i$ is a divisor with simple
  normal crossings. Write $E^{ij}= E^i\cap E^j$ and $E_{n} = \coprod
  E^{i_0\ldots i_n}$.  Let $S_0\subset X$ be set of singular points,
  $S_1\subseteq S_0$ be the set of images of $\cup E^{ij}$ and so
  on. Then the simplicial resolution is given by
$$
\xymatrix{ \ldots E_1 \sqcup S_2\ar[r]\ar@<1ex>[r]\ar@<-1ex>[r] & E_0
  \sqcup S_1\ar@<1ex>[r]\ar@<-1ex>[r] & \tilde X\sqcup S_0 }
$$
where the face maps are given by inclusions $S_i\to S_{i-1}$ on the
second component. On the first component $\delta_k$ is given by
$$
\begin{cases}
  E^{i_1\ldots i_n}\subset E^{i_i\ldots\hat i_k\ldots i_n} &\text{if $k\le n$}\\
  f: E^{i_1\ldots i_n}\to S_{n-1} &\text{if $k=n+1$}
\end{cases}
$$
\end{ex}
 
Given a simplicial resolution,
the spectral sequence (\ref{eq:ss1}) will then converge to
$H^*(X,A)$. More generally for any sheaf, there is an isomorphism
$$H^i(X,\F) \cong H^i(X_\dt,\pi_\dt^*\F)
$$
for any sheaf $\F$ on $X$.  The last property goes by the name of
{\em cohomological descent}.

Given a closed subvariety $\iota:Z\subset X$, there exists
 simplicial resolutions $Z_\dt \to Z$, $X_\dt\to X$ and a morphism $\iota_\dt:Z_\dt\to X_\dt$ covering  $\iota$. Then there is a new smooth simplicial variety $cone(\iota_\dt)$ (\cite[\S 6.3]{deligne}, \cite[IV \S 1.7]{gnpp}) whose geometric realization is 
homotopy equivalent to $X/Z$. So that the spectral sequence
\eqref{eq:ss1} converges to $H^*_c(X-Z,A)$.
 Although simplicial resolutions are far from unique, the
filtration on $H_c^*(X-Z,A)$ is the weight filtration $W$ \cite{gs},
 and this is canonically determined by $X-Z$ alone. When $A=\Q$,
this part of the datum of the canonical  mixed structure.

Let $X$ be a proper variety with a possibly empty closed set $Z$.
Let $U= X-Z$. Choose a  simplicial
resolution $C_\dt = Cone(Z_\dt\to X_\dt)$ as above.  
By convention $W$ is an increasing filtration
indexed so that
$$ W_qH_c^{p+q}(U)/W_{q-1} = E_\infty^{pq}\> (\cong E_2^{pq}\text{ over }\Q)$$
In particular, $W_{-1}=0$. The part of interest $W_0$,
can be computed as follows.  We can
form a simplicial set  by applying the connected component functor
$\pi_0$ to $C_\dt$.  This simplicial set $|\pi_0(C_\dt)|$
is called the dual complex or nerve of the simplicial resolution.
We have
$$W_0H_c^i(U,\Q)  \cong H^i(\ldots\to H^0(C_p,\Q)\to H^0(C_{p+1},\Q)\ldots )\cong H^i(|\pi_0(C_\dt)|,\Q)$$
For integer coefficients, $W_0H^i(U,\Z) = \pi^*H^i(|\pi_0(C_\dt)|,\Z)$
where $\pi:C_\dt\to \pi_0(C_\dt)$ is the constant map on components.
So that this piece of the filtration is determined by the underlying
combinatorial information encoded by the dual complex.

\begin{thm}
There is a collection of subspaces   $KH_c^i(X)\subseteq H_c^i(X,\Z)$ 
satisfying axioms (K1)-(K3) given in the introduction. Moreover, it 
is uniquely characterized by axioms.
\end{thm}

\begin{proof}
For existence, we note that $KH_c^i(X) = W_0H^i_c(X)$ satisfies these axioms
by \cite[\S 3.1]{gs}. (For rational coefficients, this goes back to 
 \cite{deligne}.)

So it remains to check uniqueness. We already checked this when $X$ is complete. The nonsingular case follows from this.
If $X$ is nonsingular, we can choose a nonsingular compactification $\bar X$.
Then from the axioms, we get 
$KH_c^i(X) = im[KH^{i-1}(\bar X-X)]$.
Then the general case now follows from the main theorem of \cite{gn} together
with remark~\ref{rmk:mayer}.

\end{proof}

From the formula $K=W_0$, we can deduce further properties.

\begin{cor}
$KH_c^i(X\times Y) = \bigoplus_{j+k=i} KH_c^j(X)\otimes KH_c^k(Y)$
\end{cor}

\begin{proof}
  \cite[thm 3]{gs}.
\end{proof}

\begin{cor}
  Let $\pi:\tilde X\to X$ be a resolution of a complete variety
  such that the exceptional divisor $E$ has normal crossings. Let $S=\pi(E)\subset X$. Then
  $\dim KH^i(X)$ is the $(i-1)$st Betti number $b_{i-1}$ of the dual
  complex of $E$  when $i> 2\dim(S)+1$.  If $S$ is
  nonsingular, then this holds for $i>1$.  When $i=2\dim(S)+1$, $\dim
  KH^i(X) = b_{i-1}$ minus the number of irreducible components of
  $S$ of maximum dimension.
\end{cor}

\begin{proof}
  This follows from lemma~\ref{lemma:mayer}, the identification of
$KH^i(E)=W_0H^i(E)$ and the above remarks. 
\end{proof}

When $X$ is a divisor with simple normal crossings, $KH_c^i(X)$ is the
cohomology of the dual complex.  As
remarked earlier \ref{rmk:ncd}, we can use a construction to a build a simplicial
variety canonically attached to $X$, for any $X$.  If we apply
$\pi_0$ to this simplicial variety, we get a simplicial set $\Sigma_X$
canonically attached to $X$, that we will call the nerve or dual
complex.  There is a canonical map $H^i(|\Sigma_X|,\Q)\to H^i(X,\Q)$
coming from the spectral sequence \eqref{eq:ss1} associated to this
simplicial variety. From the discussion in \ref{ex:ncd} and
\ref{rmk:ncd}, we can see that:

\begin{lemma}
  If $X$ is complete, the image $H^i(|\Sigma_X|,\Q)\to H^i(X,\Q)$ lies in $KH^i(X,\Q)$.
  If $X$ satisfies the assumptions of example \ref{ex:ncd}, then these
  subspaces coincide.
\end{lemma}

\section{Bounds on the combinatorial part}

Suppose that $X$ is a complete variety. Then in addition to the weight
filtration
$H^i(X,\C)$ carries a second filtration, called the Hodge filtration
induced on the abutment $H^i(X,\Omega_{X_\dt}^\dt)\cong H^i(X,\C)$ of
the spectral sequence (\ref{eq:ss3}) for $\Omega_{X_\dt}^\dt$.  By
convention $F$ is decreasing. We have $F^0= H^i(X,\C)$ and
$$F^0H^i(X,\C)/F^1 \cong H^i(X_\dt,\OO_{X_\dt})$$
$W$ induces the same filtration on the right as the one coming from (\ref{eq:ss2}). In
particular,
$$W_0 Gr^0_FH^i(X,\C) =  H^i(\ldots\to H^0(X_p,\OO)\to H^0(X_{p+1},\OO)\ldots )$$
$$\cong H^i(|\pi_0(X_\dt)|,\C) \cong W_0H^i(X,\C)$$
This means that Hodge filtrations becomes trivial on $W_0H^i(X)$. So
that this is a vector space and nothing more.

\begin{thm}
  \begin{enumerate}
  \item[]
  \item[(a)] If $X$ is a complete variety, then there is an inclusion
    $KH^i(X,\C)\hookrightarrow H^i(X,\OO_X)$.
  \item[(b)] If $f:X\to Y$ a proper morphism of varieties, then there is an
    inclusion $KH^i(f^{-1}(y),\C)\hookrightarrow
    (R^if_*\OO_X)_y\otimes \OO_y/m_y$ for each $y\in Y$.
  \end{enumerate}
\end{thm}

\begin{proof}
  The canonical map $\kappa$ factors
$$
\xymatrix{
  H^i(X,\C)\ar^{\kappa}[rr]\ar[rd] &  & H^i(X,\OO_X) \\
  & Gr^0_FH^i(X,\C)\ar[ru] & }
$$
Thus
$$W_0H^i(X,\C)\subseteq Gr^0_FH^i(X,\C)= im[ H^i(X,\C)\to  H^i(X,\OO_X)]$$
which implies (a).

For (b), let $X_y$ be the reduced fibre over $y$, and $X_y^{(n)}$ the fibre
with its $n$th infinitesimal structure. From (a), we have a
natural inclusion $s:W_0H^i(X_y,\C)\hookrightarrow
H^i(X_y,\OO_{X_y})$.  After choosing a simplicial resolution of the fibre
$f_\dt:\X_{\dt}\to X_y$, $s$ can be identified with the composition
$$E^{i0}_2(\C)\to E^{i0}_2(\OO_{\X_\dt})\to H^i(X_y,\OO_{X_y})$$
where the first map is induced by the natural map $\C\to \OO$, and the last
map is the edge homomorphism.  Applying the same construction to
the simplicial sheaf $f_\dt^*\OO_{X_y^{(n)}}$ yields a map $s_n$
fitting into a commutative diagram
$$
\xymatrix{
  W_0H^i(X,\C)\ar^{s}[rr]\ar^{s_n}[rd] &  & H^i(X_y,\OO_{X_y})\\
  & H^i(X_y,\OO_{X_y^{(n)}})\ar[ru] & }
$$
Furthermore, these maps are compatible, thus they pass to map $s_\infty$ to the
limit. Together with the formal functions theorem \cite[III
11.1]{hartshorne}, this yields a commutative diagram
$$
\xymatrix{
  W_0H^i(X,\C)\ar^{s}[rr]\ar^{s_\infty}[d]\ar^{s'}[rrd] &  & H^i(X_y,\OO_{X_y})\\
  \varprojlim H^i(X_y,\OO_{X_y^{(n)}})\ar^{\sim}[r] &
  (R^if_*\OO_X)_y\hat{}\ar[r] & (R^if_*\OO_X)_y\otimes \OO_y/m_y\ar[u]
}
$$
Since $s$ is injective, the map labeled $s'$ is injective as well.
\end{proof}

\begin{rmk}
  In item (a), we actually proved the sharper statement
$$W_0H^i(X,\C)\hookrightarrow 
Gr^0_FH^i(X,\C) = im[ H^i(X,\C)\to H^i(X,\OO_X)]$$ For certain classes
of singularities called Du Bois singularities \cite[\S 7.3.3]{ps}, which include
rational singularities \cite{kovacs}, $Gr^0_FH^i(X,\C) =
H^i(X,\OO_X)$. But this is not true in general.
\end{rmk}

\begin{cor}
  Suppose that $f:X\to Y$ is a resolution of singularities.
  \begin{enumerate}
  \item If $Y$ has rational singularities then
    $W_0H^i(f^{-1}(y),\C)=0$ for $i>0$.
  \item If $Y$ has isolated normal Cohen-Macaulay singularities,
    $W_0H^i(f^{-1}(y),\C)=0$ for $ 0< i < \dim Y-1$
  \end{enumerate}

\end{cor}

\begin{proof}
  The first statement is an immediate consequence of the theorem.  The
  second follows from the well known fact given below. We sketch the
  proof for lack of a suitable reference.

\begin{prop}
  If $f:X\to Y$ is a resolution of a variety with isolated normal
  Cohen-Macaulay singularities, then $R^if_*\OO_X = 0$ for $0<i<\dim
  Y-1$
\end{prop}

\begin{proof}[Sketch]
  We can assume that $Y$ is projective. By the Kawamata-Viehweg
  vanishing theorem \cite{kawamata,viehweg}
  \begin{equation}
    \label{eq:kV}
    H^i(X, f^*L^{-1})=0,\quad i < \dim Y = n, 
  \end{equation}
  where $L$ is ample.  Replace $L$ by $L^N$, with $N\gg 0$.  Then by
  Serre vanishing and Serre duality (we use the CM hypothesis here)
  \begin{equation}
    \label{eq:SV}
    H^i(Y, L^{-1}) = H^{n-i}(Y, \omega_Y\otimes L) = 0,\quad i < n.  
  \end{equation}
  The Leray spectral sequence together with~\eqref{eq:kV} and
  \eqref{eq:SV} imply
$$ H^0(R^if_*\OO_X\otimes L^{-1}) = 0,\quad i< n-1$$
Since the sheaves $R^if_*\OO_X$ have zero dimensional support, the
proposition follows.
\end{proof}

\end{proof}

\end{document}